\documentclass[12pt,reqno]{article}

\usepackage[usenames]{color}
\usepackage{amssymb}
\usepackage{amsmath}
\usepackage{amsthm}
\usepackage{amsfonts}
\usepackage{amscd}
\usepackage{graphicx}

\usepackage[colorlinks=true,
linkcolor=webgreen,
filecolor=webbrown,
citecolor=webgreen]{hyperref}

\definecolor{webgreen}{rgb}{0,.5,0}
\definecolor{webbrown}{rgb}{.6,0,0}
\usepackage{color}
\usepackage{fullpage}
\usepackage{float}

\usepackage{cases}

\usepackage{graphics,amsmath,amssymb}
\usepackage{amsthm}
\usepackage{amsfonts}
\usepackage{latexsym}
\usepackage{enumitem}

\begin{document}

\theoremstyle{plain}
\newtheorem{theorem}{Theorem}[section]
\newtheorem{remark}{Remark}[section]
\newtheorem{proposition}{Proposition}[section]
\newtheorem{example}[remark]{Example}
\newtheorem{lemma}{Lemma}
\newtheorem{definition}{Definition}
\newtheorem{corollary}{Corollary}

\renewcommand{\thefootnote}{\fnsymbol{footnote}}

\begin{center}
\vskip 1cm{\Large\bf The most concise recurrence formula for the sums
\vskip .08in of integer powers}\footnote[5]{
This paper has been originally published in {\it The Mathematical Gazette} 109 (July 2025) \mbox{pp.\ 253--258;} \url{https://doi.org/10.1017/mag.2025.10073}. In this arXiv version, we provide an abstract and correct a few minor misprints in the published version.}
\vskip .3in \large Jos\'{e} Luis Cereceda \\
{\normalsize Collado Villalba, 28400 (Madrid), Spain} \\
\href{mailto:jl.cereceda@movistar.es}{\normalsize{\tt jl.cereceda@movistar.es}}
\end{center}

\begin{abstract}
For integers $n,k \geq 1$, let $S_k(n)$ denote the power sum $1^k +2^k + \cdots + n^k$. In this note, we first recall the minimal recurrence relation connecting $S_k(n)$ and $S_{k-1}(n)$ established by Abramovich (1973). We then discuss an old algorithm to determine the coefficients of the power sum polynomial $S_k(n)$ in terms of the coefficients of $S_{k-1}(n)$ (see, e.g., Bloom (1993) and Owens (1992)). Moreover, we bring to light an explicit relationship between $S_k(n)$ and $S_{k+1}(n)$ put forward by Budin and Cantor (1972). We conclude that these procedures (including the integration formula expressing $S_k(n)$ in terms of $S_{k-1}(n)$) all constitute equivalent methods to determine $S_k(n)$ starting from $S_{k-1}(n)$. In addition, as a by-product, we provide a determinantal formula for the Bernoulli numbers involving the binomial coefficients.
\end{abstract}

\vspace{.5cm}

For integers $n,k \geq 1$, let $S_k(n)$ denote the power sum $1^k + 2^k + \cdots + n^k$, with $S_0(n) =n$. As is well known, $S_k(n)$ can be represented by a polynomial in $n$ of degree $k+1$ with constant term equal to zero, i.e. $S_k(n) = \sum_{j=1}^{k+1} a_{k,j} n^j$, for certain rational coefficients $a_{k,j}$. In 1973, in the popular science magazine {\it Kvant\/}\footnote{
This magazine was founded in 1970 in the former Soviet Union. The prominent mathematician \mbox{A. N. Kolmogorov} (1903-1987) played an active role in the inception of the magazine, and published various articles in it.}, appeared an article by Vladimir Abramovich \cite{abra} in which, among other things, the author derived the following minimal recurrence relation which expresses in a very compact way the connection between the power sums $S_k(n)$ and $S_{k-1}(n)$:
\begin{equation}\label{abra}
S_k(n) = n + k S_{k-1}^{\star}(n), \quad n,k \geq 1,
\end{equation}
where the polynomial $S_{k-1}^{\star}(n)$ is constructed by replacing in $S_{k-1}(n)$ the term $n^j$ with the expression $\frac{n^{j+1}-n}{j+1}$ for each $j =1,2,\ldots,k$. Note that $S_{k-1}^{\star}(n)$ has degree $k+1$, as it should be.

We can use \eqref{abra} to get $S_k(n)$ starting from $S_{k-1}(n)$. As a simple example, using the fact that $S_1(n) = \frac{1}{2}n^2 + \frac{1}{2}n$, we can find $S_2(n)$ as follows:
\begin{equation*}
S_2(n) = n +2 S_1^{\star}(n) = n + \frac{n^3 -n}{3} + \frac{n^2 -n}{2} = \frac{n(n+1)(2n+1)}{6}.
\end{equation*}

The proof of \eqref{abra} given in \cite{abra} employed a rather complicated mathematical induction. A different proof of \eqref{abra} was subsequently given by the same author\footnote{The mathematician Vladimir Samuel Shevelev (1945-2018), originally named Vladimir Abramovich, adopted the surname Shevelev in 1982.} in \cite{shevelev} with the help of the Bernoulli polynomials. In this note we point out that the recurrence formula \eqref{abra} is actually equivalent to another concise recurrence formula for $S_{k}(n)$, namely the formula in \eqref{sher} below. Furthermore, we show that the procedure to get $S_k(n)$ from $S_{k-1}(n)$ entailed by \eqref{abra} amounts to a seemingly unrelated method to determine the coefficients of the power sum polynomial $S_k(n)$ in terms of the coefficients of $S_{k-1}(n)$. This method is discussed in, e.g. Bloom \cite{bloom} and Owens \cite{owens}, and in this note we shall refer to it as the Bloom-Owens algorithm. Finally, as a by-product, we provide a determinantal formula for the Bernoulli numbers involving the binomial coefficients (see \eqref{det} below).

\section{Proof of recurrence \eqref{abra}}

Let $S_k(x)$ denote the polynomial obtained by replacing $n$ by $x$ in $S_k(n)$, i.e. $S_k(x) = \sum_{j=1}^{k+1} a_{k,j} x^j$, where $x$ is a real variable. To derive \eqref{abra}, we make use of the following integration formula for $S_k(x)$:
\begin{equation}\label{sher}
S_k(x) = k \int_{0}^{x} S_{k-1}(t) dt \, + \left( 1 - k \int_{0}^{1} S_{k-1}(t) dt \right) x, \quad k \geq 1,
\end{equation}
with $S_0(x) = x$. The formula above constitutes a well-established, classical result, and can be taken as the definition of $S_k(x)$. (For a proof of \eqref{sher} see, e.g., \cite{sherwood}, Theorem 3.4.1 of \cite{coen}, and \cite{costa}. In the latter reference \cite{costa}, a geometrical interpretation of \eqref{sher} is also given.) Moreover, the integration formula \eqref{sher} determines $S_k(x)$ entirely in terms of $S_{k-1}(x)$, making it suitable for comparison with the recurrence formula \eqref{abra}. Indeed, it is immediate to check that \eqref{sher} can be cast into the form $S_k(x) = x + k S_{k-1}^{\star}(x)$ provided that
\begin{equation*}
S_{k-1}^{\star}(x) = \int_{0}^{x} S_{k-1}(t) dt - x  \int_{0}^{1} S_{k-1}(t) dt.
\end{equation*}
Hence, substituting $S_{k-1}(t)$ by the polynomial $\sum_{j=1}^{k} a_{k-1,j} t^j$, it follows that
\begin{equation*}
S_{k-1}^{\star}(x) = \sum_{j=1}^k a_{k-1,j} \left( \int_{0}^x t^j dt - x \int_{0}^1 t^j dt \right)
= \sum_{j=1}^k a_{k-1,j} \frac{x^{j+1} -x}{j+1}.
\end{equation*}
Thus, in the case that $x =n$ is a natural number, we retrieve \eqref{abra} and the proof is complete. Conversely, as can be easily shown, the integration formula \eqref{sher} can be deduced from Abramovich's formula \eqref{abra}.

\section{The Bloom-Owens algorithm}

There exists a simple procedure to generate the coefficients $a_{k,j}$ of the power sum polynomial $S_k(n)$ once the coefficients $a_{k-1,j}$ of $S_{k-1}(n)$ are known. The method (which has been nicely described by Bloom \cite{bloom} and Owens \cite{owens}; see also \cite{carchidi,doucet,torabi}) consists of the following three steps:
\begin{enumerate}[label=(\roman*)]
  \item Start with $S_{k-1}(n) = \sum_{j=1}^k a_{k-1,j} n^j$ with known coefficients $a_{k-1,j}$, $j=1,2,\ldots,k$.
  \item For each $j=2,3,\ldots,k+1$, set $a_{k,j} = \frac{k}{j} \, a_{k-1,j-1}$. In this way, all the coefficients of $S_k(n)$ except for $a_{k,1}$ are determined.
  \item Since $S_k(1) =1$ the coefficient $a_{k,1}$ is obtained from the equation $a_{k,1} = 1- \sum_{j=2}^{k+1} a_{k,j}$, where the coefficients $a_{k,j}$ (for $j =2,3,\ldots,k+1$) are those determined in step (ii).
\end{enumerate}

We can confirm that, as one might expect, the mathematical operation behind the recurrence formula \eqref{abra} is precisely the same as that implemented by the Bloom-Owens algorithm embodied in steps (i), (ii), and (iii) listed above. Indeed, using the definition of $S_{k-1}^{\star}(n)$, from \eqref{abra} we obtain
\begin{align*}
\sum_{j=1}^{k+1} a_{k,j} n^j & = n + k \, \sum_{j=1}^k a_{k-1,j} \frac{n^{j+1} -n}{j+1} \\
& = k \, \sum_{j=2}^{k+1} \frac{a_{k-1,j-1}}{j} n^j + \Bigg( 1 -k \sum_{j=2}^{k+1} \frac{a_{k-1,j-1}}{j}\Bigg)n.
\end{align*}
Now, comparing the coefficients of $n^j$ in the left and right sides of the preceding equation, we find that
\begin{equation*}
a_{k,j} = \frac{k}{j} \, a_{k-1,j-1} \quad\text{for}\,\,  j =2,3,\ldots,k+1,
\end{equation*}
and
\begin{equation}\label{boa}
a_{k,1} = 1 - k \sum_{j=2}^{k+1} \frac{a_{k-1,j-1}}{j} = 1 - \sum_{j=2}^{k+1} a_{k,j},
\end{equation}
yielding the same output as that prescribed in steps (ii) and (iii). We therefore conclude that Abramovich's formula \eqref{abra} and the Bloom-Owens algorithm are, in fact, equivalent.

Incidentally, as pointed out by Bloom \cite{bloom}, the method just described is treated somewhat peripherally in Edwards' book \cite{edw}. It should be added, however, that Edwards' account of the method is done in terms of the Bernoulli polynomials, and there is no explicit mention of steps (i), (ii), and (iii) given above. Perhaps of more interest to us is the following proposition regarding the power sum polynomials $S_k(n)$ and $S_{k+1}(n)$ put forward by Budin and Cantor in a short note of 1972 in the IEEE journal \emph{Transactions on Systems, Man, and Cybernetics} \cite{budin} (in our notation):

\noindent
{\it Proposition} (Budin \& Cantor, 1972): If
\begin{equation*}
S_k(n) = b_1 n^{k+1} + b_2 n^k + b_3 n^{k-1} + \cdots + b_{k+1}n,
\end{equation*}
then
\begin{equation}\label{budin}
\begin{split}
S_{k+1}(n) & = \frac{k+1}{k+2}b_1 n^{k+2} + \frac{k+1}{k+1}b_2 n^{k+1} + \frac{k+1}{k}b_3 n^{k} + \cdots \\
& \qquad + \frac{k+1}{2}b_{k+1} n^{2} + \bigg[ 1 - (k+1) \sum_{i=1}^{k+1} \frac{b_i}{k+3-i} \bigg] n, \quad k \geq 0.
\end{split}
\end{equation}

Formula \eqref{budin} allows one to compute $S_{k+1}(n)$ starting from $S_k(n)$ (or, by the way, $S_{k}(n)$ starting from $S_{k-1}(n)$). As can be easily verified, the involved procedure corresponds exactly with the Bloom-Owens algorithm.

In conclusion, we have shown that the Abramovich formula \eqref{abra}, the integration formula \eqref{sher}, the Bloom-Owens algorithm, and the Budin-Cantor formula \eqref{budin} all constitute equivalent methods for obtaining the coefficients of $S_{k}(n)$ in terms of the coefficients of $S_{k-1}(n)$.

\section{Concluding remarks}

To finish, it is worth recalling the well-known fact that the coefficients $a_{k,j}$ in $S_k(n) = \sum_{j=1}^{k+1} a_{k,j} n^j$ are intimately related to the Bernoulli numbers. This is best appreciated when we write $S_k(n)$ in the Bernoulli form (see, e.g. \cite{wu}),
\begin{equation*}
S_k(n) = \frac{1}{k+1} \sum_{j=1}^{k+1} (-1)^{k+1-j} \binom{k+1}{j} B_{k+1-j} n^j,
\end{equation*}
where $B_0 =1$, $B_1 = -\frac{1}{2}$, $B_2 =\frac{1}{6}$, $B_3 =0$, $B_4 = -\frac{1}{30}$, etc., are the Bernoulli numbers. The coefficients $a_{k,j}$ are then given explicitly by
\begin{equation*}
a_{k,j} = \frac{(-1)^{k+1-j}}{k+1} \binom{k+1}{j} B_{k+1-j}, \quad j =1,2,\ldots,k+1.
\end{equation*}
Note, in particular, that $a_{k,1} = (-1)^k B_k$. In view of \eqref{sher}, this implies that
\begin{equation*}
(-1)^k B_k = 1 - k \int_{0}^{1} S_{k-1}(t) dt,  \quad k \geq 1.
\end{equation*}

On the other hand, upon substitution of the above explicit expression of $a_{k,j}$ into \eqref{boa}, we quickly arrive at the following recurrence relation for the Bernoulli numbers:
\begin{equation}\label{recb}
\sum_{j=1}^{k} (-1)^j \binom{k+1}{j} B_j = k, \quad k \geq 1,
\end{equation}
yielding $B_1 = -\frac{1}{2}$ for $k=1$. When $k \geq 2$, we can use \eqref{recb} to get $B_k$ in terms of $B_1, B_2, \ldots, B_{k-1}$.

Lastly, considering \eqref{recb} as a system of $k$ linear equations in the variables $B_1, B_2, \ldots, B_k$, and applying Cramer's rule to such a system, we are able to deduce the following determinantal formula for $B_k$:
\begin{equation}\label{det}
B_{k} = \frac{(-1)^k}{(k+1)!} \!
\begin{vmatrix}
\binom{2}{1} & 0 & 0 & \! \hdots & 0 & \! 1 \\[3pt]
\binom{3}{1} & \! \binom{3}{2} & 0 & \! \hdots & 0 & \! 2 \\[3pt]
\binom{4}{1} & \! \binom{4}{2} & \! \binom{4}{3} & \! \hdots  & 0 & \! 3 \\[3pt]
\vdots & \vdots & \vdots & \! \ddots & \vdots & \vdots \\[3pt]
\binom{k}{1} & \! \binom{k}{2} & \! \binom{k}{3} & \!\hdots & \! \binom{k}{k-1} &
\! k-1 \\[5pt]
\binom{k+1}{1} & \! \binom{k+1}{2} & \! \binom{k+1}{3} & \! \hdots & \! \binom{k+1}{k-1} &
\! k
\end{vmatrix},
\end{equation}
with the convention that $B_1 = -\frac{1}{2}$. Furthermore, by using the basic properties of determinants, the above formula can be transformed so as to obtain the following determinantal formula of order $k-1$ for $B_k$:
\begin{equation*}
B_{k} = \frac{1}{(k+1)!} \!
\begin{vmatrix}
1 & \! \binom{3}{2} & \! 0 & \! \hdots & 0 & \! 0 \\[3pt]
2 & \! \binom{4}{2} & \! \binom{4}{3} & \! \hdots  & 0 & \! 0 \\[3pt]
3 & \! \binom{5}{2} & \! \binom{5}{3} & \! \hdots  & 0 & \! 0 \\[3pt]
\vdots & \vdots & \vdots & \! \ddots & \vdots & \vdots \\[3pt]
k-2 & \! \binom{k}{2} & \! \binom{k}{3} & \!\hdots & \! \binom{k}{k-2} &
\! \binom{k}{k-1} \\[5pt]
k-1 & \! \binom{k+1}{2} & \! \binom{k+1}{3} & \! \hdots & \! \binom{k+1}{k-2} &
\! \binom{k+1}{k-1}
\end{vmatrix}, \quad k \geq 2.
\end{equation*}
For example, for $k=8$, this last formula gives us
\begin{equation*}
B_{8} = \frac{1}{9!} \,
\begin{vmatrix}
1 & 3 & 0 & 0 & 0 & 0 & 0 \\
2 & 6 & 4 & 0 & 0 & 0 & 0 \\
3 & 10 & 10 & 5 & 0 & 0 & 0 \\
4 & 15 & 20 & 15 & 6 & 0 & 0 \\
5 & 21 & 35 & 35 & 21 & 7 & 0 \\
6 & 28 & 56 & 70 & 56 & 28 & 8 \\
7 & 36 & 84 & 126 & 126 & 84 & 36 \\
\end{vmatrix}
= -\frac{1}{30}.
\end{equation*}

Finally, we may remark, that a determinantal formula for $B_k$ equivalent to \eqref{det} can be found (along with a historical comment) in \cite{booth}. For other determinantal expressions of the Bernoulli numbers and polynomials, the reader may consult \cite{qi} and references therein.

\vspace{.3cm}

\noindent \emph{Acknowledgements}

{\small
The author would like to thank the anonymous reviewer for useful comments and suggestions that led to improvements in an earlier version of this article.}

\end{document}